\documentclass{article}

\usepackage{amsfonts}

\usepackage{amssymb}

\usepackage{graphicx}

\newtheorem{theorem}{Theorem}
\newtheorem{proposition}[theorem]{Proposition}

\newtheorem{definition}[theorem]{Definition}

\begin{document}

\title{Nobody needs equations}
\author{Adonai S. Sant'Anna\\ {\small adonaisantanna@gmail.com}}  

\maketitle

\begin{abstract}
Equations are ubiquitous in most mathematical activities. Nevertheless, in this paper it is shown how to do standard mathematics without any equation at all. More than that, it is proven there is a foundational framework for standard mathematics where equations cannot be even written, in the sense they are not formulas. The proof of those claims is very simple, almost obvious. I use this framework to suggest a way to deal with certain notions of indiscernibility between `objects', with special emphasis on some aspects of quantum mechanics. Finally I compare this approach to quasi-set theory, an unnecessarily complicated formal work designed to deal with violation of Leibniz’s Principle of the Identity of Indiscernibles.
\end{abstract}

\section{Introduction}

The main claim here is this: most of standard mathematics can be done without the need or even the possibility of writing any equation at all. Before any proof of that claim is provided, a couple of concepts needs clarification: What is standard mathematics? What is an equation?

Due to a vast plurality of mathematical ideas, concepts, and theories, there is no clear-cut definition of what standard mathematics is supposed to be. A simple way to illustrate such a diversity of mathematical ideas is the Mathematics Subject Classification produced by the editorial staff of Mathematical Reviews and Zentralblatt MATH. In that document we always find, in each discipline of mathematics, a subsection that reads ``None of the above, but in this section''. Thus, mathematics is, in principle, open to new possibilities. So, I appeal to a sociological notion which is supposed to reflect the basic everyday mathematics. Within this scope, standard mathematics refers to those branches of knowledge where a significant part of the mathematical community works, including undergraduate and graduate studies focused on mathematics, pure or not. Well known examples are algebra, topology, real and complex analysis, differential equations, numerical analysis, probability theory, measure theory, functional analysis, linear algebra, distribution theory, combinatorics, number theory, logic and foundations, multivariate analysis, mathematical methods in physics, and so on. With this notion in hand the reader could easily say ``well, standard mathematics employs equations; therefore, this paper's abstract is inconsistent''. To answer such a possible criticism, next paragraphs are essential.

Contrary to the informal notion of standard mathematics, the concept of {\em equation\/} can be given in a quite clear fashion. An equation is an atomic formula $x = y$, where $x$ and $y$ are terms, and $=$ is a very specific binary predicate. In other words, we are talking about a formal language, where equations are particular cases of formulas. The most usual formal language for that purpose is that one of Zermelo-Fraenkel set theory (and its variations like ZFC, ZFA, ZF $+$ CH, ZF $+$ inaccessible cardinals, etc.) \cite{Jech-03, Button-18}.

ZF has the well known advantage of being powerful enough to provide foundations for different branches of mathematics, like the ones (most of them) already listed above. So, if someone shows a variation of ZF with no identity, somehow equivalent to ZF, then that entails all those branches of mathematics can be done without any equation at all. That's the strategy to answer how to do standard mathematics with no equations. So, the first answer to the possible criticism above is this: it can be provided a formal framework as a foundation for standard mathematics where this new framework is equivalent (in a precise sense) to ZF. Besides, within this new framework there is no equation at all.

Obviously it is not possible to simply eliminate identity from ZF and expect everything works nicely from there. In order to cope with that difficulty it is necessary to select appropriate models of our new version of ZF and let them do the work. Thus, some basic stuff on model theory is required in this paper. Otherwise, there is no simple way to get rid of equations.

ZF is a first-order theory {\em with\/} identity \cite{Mendelson-97}. In Mendelson's book \cite{Mendelson-97} the reader can find precise definitions of formal languages, formal theories, first-order languages, first-order theories, first-order theories with identity, terms, and formulas, although the formal set theory examined in that book is NBG. A nice reference for a detailed review on ZF and some of its variations is Jech's book \cite{Jech-03}. Once this formal language is settled, we can finally state in a more detailed way: an equation is a formula

$$x = y,$$

\noindent where $x$ and $y$ are terms, and $=$ is a binary predicate where reflexivity, transitivity, symmetry and substitutivity hold. A term, on its turn, can be either a constant, or a variable or a functional letter applied to terms. The definition of terms and formulas is given recursively. To {\em solve\/} an equation is to determine which possible values of $x$ and $y$ grant the equation $x=y$ to be a theorem, where the concept of theorem is given as well in Mendelson's book. Eventually an equation has only one solution, several solutions or no solution.

So, following this brief introduction, the next step is to recall the proper axioms of ZF. Since we need some basic results from model theory \cite{Jech-03}, a brief review about that is provided as well. The last step is to show how a different version of ZF (without identity) can be used as a possible foundation for standard mathematics.

All this very simple mathematical apparatus is finally used to show that quasi-set theory \cite{French-06} is an unnecessarily complicated way to cope with the problem of non-individuality, in the sense of a violation of Leibniz’s Principle of the Identity of Indiscernibles (LPII). The idea of violation of that principle has been contemplated by many philosophers since the dawn of quantum mechanics, a theory which seems to demand the existence of indiscernible objects that are not necessarily identical. That is an important issue, since some authors still keep insisting over the last decades that (at least) certain aspects of quantum mechanics demand some sort of non-standard mathematics to violate LPII \cite{Krause-92}\cite{Arenhart-12}\cite{Barros-23}\cite{Domenach-08}\cite{French-06}\cite{Chiara-93}. In this paper I show a rather simple modification of ZF, working together with model theory, can handle the most common philosophical motivations usually employed to advocate quasi-set theory. A non-standard approach like quasi-set theory \cite{French-06} is excessively complicated to deal with violation of LPII. Besides, it is not mathematically clear that that framework does the job \cite{Sant'Anna-23}. It certainly does not provide mathematical foundations for certain physical phenomena \cite{Sant'Anna-20}.

In a nutshell, our proposal is this: 

\begin{enumerate}

\item A weaker version of ZF, called $Z\Phi$, is introduced, without identity, without the Axiom of Extensionality;

\item Two metatheorems grant that, for any transitive model of $Z\Phi$, there is a corresponding model of ZF; and, conversely, any model of ZF is a transitive model of $Z\Phi$; thus, transitive models of $Z\Phi$ provide a foundation for standard mathematics;

\item Instead of using $Z\Phi$ as a foundation for mathematics, we use only its transitive models for that purpose; those models do not need to be endowed with identity, in the sense of a translation of the binary predicate $=$ in first order languages; after all, transitivity can be expressed by means of membership alone; henceforth, $Z\Phi$ is used here as a starting point to introduce certain models which, in their turn, can be used as a foundational tool;

\item If we assume that quantum mechanics demands standard mathematics and, besides that, a violation of LPII at least for some objects, then we use any model from a specific class of non-transitive models of $Z\Phi$ who are capable to reconcile both assumptions;
    
\item Last item above is accomplished thanks to the fact that some non-transitive models of $Z\Phi$ admit submodels that are transitive;
    
\item $Z\Phi$ demands less concepts and less axioms than quasi-set theory;

\item We sketch a proof that quasi-set theory might provide a way to deal with standard mathematics and still violate LPII for some objects.

\end{enumerate}

\section{ZF}\label{zf}

ZF is a first-order theory with identity ($=$) and a binary predicate letter $\in$ where we read $x\in y$ as `$x$ belongs to $y$' (or, alternatively, `$x$ is a member of $y$'). In this sense, ZF has two primitive concepts, namely, identity $=$ and membership $\in$. The negation of $x=y$ and $x\in y$ are denoted by $x\neq y$ and $x\not\in y$, respectively. Besides its logical axioms (the same as NBG in Mendelson's book), the proper axioms of ZF are as follows.

Sets who share the same elements are identical.

\begin{description}

\item[ZF1 - Extensionality] $$\forall x\forall y\forall z((z\in x \Leftrightarrow z\in y)\Rightarrow x = y).$$

\end{description}

There is a set with no element.

\begin{description}

\item[ZF2 - Empty Set] $$\exists x\forall y(y\not\in x).$$

\end{description}

For any sets $x$ and $y$ there is a set $z$ whose only elements are $x$ and $y$.

\begin{description}

\item[ZF3 - Pair] $$\forall x\forall y \exists z \forall t(t\in z \Leftrightarrow (t = x \vee t = y)).$$

\end{description}

For any set $x$ there is a set whose elements are the subsets of $x$.

\begin{description}

\item[ZF4 - Power] $$\forall x\exists y\forall t(t\in y \Leftrightarrow t\subseteq x),$$

\noindent where $t\subseteq x$ is an abbreviation of formula $\forall r(r\in t\Rightarrow r\in x)$.

\end{description}

For any set $x$ there is the arbitrary union of all sets who belong to $x$.

\begin{description}

\item[ZF5 - Union] $$\forall x\exists y\forall z(z\in y \Leftrightarrow \exists w(z\in w \wedge w\in x)).$$

\end{description}

A set $x$ can be defined by a formula as long the elements of $x$ are elements of another set $z$.

\begin{description}

\item[ZF6$_{\mathcal{F}}$ - Separation] $$\forall z\exists x\forall y(y\in x \Leftrightarrow y\in z\wedge \mathcal{F}(y)),$$

\noindent where $\mathcal{F}(y)$ is a formula with no free occurrences of $x$.

\end{description}

There is at least one infinite set.

\begin{description}

\item[ZF7 - Infinity] $$\exists x(\varnothing\in x \wedge \forall y(y\in x\Rightarrow S(y)\in x)),$$

\noindent where $\varnothing$ is a constant called `empty set' whose existence is granted by ZF2 and whose uniqueness is granted by ZF1. Besides, $S(y) = y\cup \{ y\}$, where $\cup$ is a finite union granted by ZF5, and $\{ y\}$ is a pair whose only element is $y$, granted by ZF3.

\end{description}

Certain formulas can be used to define new sets in a way quite different from the Separation Scheme.

\begin{description}

\item[ZF8$_{\mathcal{F}}$ - Replacement] $$\forall x\exists!y\mathcal{F}(x,y)\Rightarrow\forall z\exists w\forall t(t\in w \Leftrightarrow\exists s(s\in z\wedge\mathcal{F}(s,t))),$$

\noindent where $\mathcal{F}(x,y)$ is a formula where all occurrences of $x$ and $y$ are free. Besides, $\exists!y\mathcal{F}(x,y)$ is an abbreviation of formula $\exists y\exists z ((\mathcal{F}(x,y)\wedge \mathcal{F}(x,z))\Rightarrow y = z)$.

\end{description}

ZF is well-founded, i.e., there is no infinite chains of membership like $x\in y$ and $y\in z$ and $z\in x$.

\begin{description}

\item[ZF9 - Regularity] $$\forall x(x\neq\varnothing\Rightarrow \exists y(y\in x \wedge x\cap y = \varnothing)),$$

\noindent where finite intersection $\cap$ can be defined from finite union in the usual way.

\end{description}

Observe there is no occurrence of identity in axioms {\bf ZF2}, {\bf ZF4}, {\bf ZF5}, and {\bf ZF6}$_{\mathcal{F}}$. I want to take advantage of that feature of ZF. Now, how do we handle with the remaining axioms if we want to get rid of identity? We answer that in the next two Sections.

\section{Models of ZF}\label{sobremodelos}

An {\em interpretation\/} for a given language $\cal L$ is an ordered pair ${\cal U} = (U, I)$, where $U$ is the universe of ${\cal U}$ and $I$ is the interpretation function which maps the symbols of $\cal L$ to appropriate relations, functions, and constants in $U$. Such an interpretation is defined by means of another language $\cal L'$ instead of $\cal L$. The interpretation function $I$ is simply a translation from $\cal L$ to $\cal L'$. If all translated axioms of a first-order theory hold in $\cal U$, we say $\cal U$ is a {\em model\/} of that theory. For details see \cite{Jech-03}.

ZF has a lot of well-known models \cite{Jech-03}. Usually model theory is employed to provide metamathematical results about independence of formulas (with respect to certain ensembles of formulas) or metatheorems concerning relative consistency between theories. Some of the most famous results of model theory are (i) the consistency of the Axiom of Choice (AC) with respect to the remaining axioms of ZFC; (ii) the consistency of the Continuum Hypothesis (CH) with respect to the remaining axioms of ZF $+$ CH; (iii) the independence of AC; and (iv) the independence of CH. The first two were proved by Kurt G\"odel, while the last two were demonstrated by Paul Cohen.

A very interesting class of partial models (interpretations where just some axioms hold) of ZFC is that one by Abian and LaMacchia \cite{Abian-78}. Those authors introduced simple partial models of ZFC where any set is translated as a sum of powers of 2. That is possible if we ignore the Axiom of Infinity. That's why Abian and LaMacchia interpretations are just partial models: The Axiom of Infinity does not hold in any of them. Thanks to that technique they prove the independence of several axioms of ZFC. In particular, the proof of independence of the Axiom of Extensionality is provided by a partial model which admits the existence of two empty sets. That was accomplished thanks to a {\em non-transitive interpretation\/} of ZFC. That's a key point in this paper. I mention this class of models because Abian and LaMacchia beautiful four-page work can be easily understood even by those with little familiarity with set theories and model theory.

An interpretation ${\cal U} = (U, I)$ of the language of ZF is {\em transitive\/} iff, for any $x$ and any $y$, `$x$ belonging to $U$ and $y$ belonging to $x$' entail `$y$ belonging to $U$', where `belonging' is the translation of $\in$ by means of function $I$. The good news is that the Axiom of Extensionality is {\em equivalent\/} to transitive interpretations of the language of ZF. In other words, Axiom of Extensionality holds in $\cal U$ iff $\cal U$ is transitive (as long we translate ZF identity as identity in the interpretation). That well known metamathematical result \cite{Jech-03} is our golden ticket to mathematics without equations, as we can see in the next Section.

But before that, one last remark here. If ZF provides a foundation for standard mathematics, the same happens to any model of ZF. The only disadvantage of using models as foundation is that they are too restrictive. ZF opens a lot of possibilities, like those where AC holds and those where AC does not hold. Many other possibilities are available in ZF, involving the Continuum Hypothesis, the existence of inaccessible cardinals, among other independent formulas.

\section{A variation of ZF with no identity}

Axiom of Extensionality {\bf ZF1} in ZF is one possible way to translate the intuitive idea of Leibniz’s Principle of the Identity of Indiscernibles. The rationale for such a claim is straightforward. Indiscernible sets are those who share the same elements (recall the word `same' here makes reference to repeated occurrences of a term in a formula) . So, {\bf ZF1} says that indiscernible sets are identical, which is a set-theoretic version of LPII. In the next pages we show it is possible to violate LPII and still keep the standard way of doing mathematics. All we need to do is to delete {\bf ZF1}, eliminate identity, and specify which models of this new theory can be used as a foundational framework.

I refer to this new version of ZF as $Z\Phi$. I use capital Greek letters Zeta and Phi as a homage to Heraclitus, the pre-Socratic philosopher who seemed to defy the concept of identity. $Z\Phi$ (phonetically equivalent to ZF) is a first order theory {\em without\/} identity, with a binary predicate letter $\in$ where we read $x\in y$ as `$x$ belongs to $y$' (or, alternatively, `$x$ is a member of $y$'). In this sense, $Z\Phi$ has one single primitive concept, namely, membership $\in$. The negation of $x\in y$ is denoted by $x\not\in y$. Besides its logical axioms (the same as NBG in Mendelson's book), the proper axioms of $Z\Phi$ are as follows: the same proper axioms of ZF, except {\bf ZF1} (Extensionality). In the remaining axioms from Section \ref{zf}, whenever occurs a formula $x = y$ we replace it by

$$\forall t(t\in x \Leftrightarrow t\in y).$$

Whenever occurs a formula $x \neq y$ we replace it by

$$\exists t((t\in x\wedge t\not\in y)\vee (t\in y\wedge t\not\in x)).$$

I dropped {\bf ZF1} because it is the only postulate in ZF which genuinely demands identity $=$. What do I mean by that? In axioms {\bf ZF2}$\sim${\bf ZF9} any occurrence of formula $x=y$ can be replaced by $\forall t(t\in x \Leftrightarrow t\in y).$ In axiom {\bf ZF1} that is not the case. Since axiom {\bf ZF1} is independent from the remaining postulates, it describes a non-trivial relation between membership and identity. Extensionality is exactly about the relationship between identity and membership. However, if we try to use $Z\Phi$ as a foundation for standard mathematics, we will have a lot of issues. For example, the finite ordinals whose existence is granted by the Axiom of Infinity will not be necessarily unique. That entails $Z\Phi$ is consistent with the existence of many ordinals $0$, many ordinals $1$, many ordinals $2$, and so forth. Why is that? Because there may be models of $Z\Phi$ who are not transitive. After all, we got rid of {\bf ZF1}. So, how do we avoid such an inconvenience? The answer is straightforward. All we have to do is to assume that the only models of $Z\Phi$ that matter are the transitive ones. But how to discriminate between those models that matter and all the rest?

My proposal here is to use a {\em semantical approach\/} for the sake of applicability of $Z\Phi$ as a foundational framework for standard mathematics. $Z\Phi$ alone is not very useful for certain purposes (like the one that demands uniqueness of each and every ordinal). On the contrary, $Z\Phi$ is a mess since it has no postulate to grant transitivity for all its models. But if we prove standard mathematics can be done on any transitive model of $Z\Phi$, then we have a quite useful framework for standard mathematics.

To make our proposal sound we need the next metatheorems:

\begin{proposition}\label{meta}

Any model of ZF is a transitive model of $Z\Phi$.

\end{proposition}

\begin{description}

\item[\sc Proof:] Let ${\cal U} = (U, I)$ be a model of ZF. Then $I$ translates $\in$ as a relation $\in_U$ in $\cal U$ and translates identity $=$ as an identity in $\cal U$. Since axioms {\bf Z2}$\sim${\bf Z9} (all of them can be written using only membership $\in$) are the same in both ZF and $Z\Phi$ (recall ZF is a stronger version of $Z\Phi$, in the sense there is one more primitive concept and one more independent axiom in ZF), then those postulates hold in $Z\Phi$. Hence, $\cal U$ is a model of $Z\Phi$. Since {\bf ZF1} holds in $\cal U$, then such a model is transitive.

\end{description}

\begin{proposition}\label{meta2}

Any transitive model of $Z\Phi$, endowed with identity, induces a model of ZF.

\end{proposition}

\begin{description}

\item[\sc Proof:] Suppose ${\cal U} = (U, I)$ is a model of $Z\Phi$. Then $I$ translates $\in$ as a relation $\in_U$ in $\cal U$. Since axioms {\bf Z2}$\sim${\bf Z9} are the same in both ZF and $Z\Phi$, then those postulates hold in ZF. According to $I$, formula $\forall t(t\in x \Leftrightarrow t\in y)$ translates in $\cal U$ as $\forall t(t\in_U x \Leftrightarrow t\in_U y)$. But $\cal U$ has identity. Besides, $\cal U$, by hypothesis, is transitive. Thus, $\forall t((t\in_U x \Leftrightarrow t\in_U y)\Rightarrow x = y)$ holds in $\cal U$. Therefore, identity in ZF can be translated as identity in $\cal U$. Hence, this new version of $\cal U$ (there is now a translation of identity) is a model of ZF.

\end{description}

As a simple example to illustrate the consequences of last propositions, consider the next differential equation:

\begin{equation}\label{1}
\frac{dy}{dx} = y,
\end{equation}

\noindent where $y$ is a real function whose domain is the set of real numbers. Within the context of ZF such an equation is an identity between sets. After all, any function is simply a set of ordered pairs. Set $y$ is an arbitrary set of ordered pairs of real numbers (as long it is a differentiable real function), and set $\frac{dy}{dx}$ is another set of ordered pairs of real numbers called `the derivative of $y$'. Thus, still within ZF, the very same differential equation can be rewritten as it follows:

\begin{equation}\label{2}
\forall t\left( t\in \frac{dy}{dx} \Leftrightarrow t\in y\right).
\end{equation}

Recall those $t$ above are ordered pairs of real numbers, at least in the case where $t$ belongs to $y$.

ZF grants last formula entails $\frac{dy}{dx} = y$, due to the Axiom of Extensionality. Substitutivity of identity grants formula \ref{1} entails formula \ref{2}. Thus, ZF guarantees formulas \ref{1} and \ref{2} are equivalent. This fact by itself is enough to prove equations are unnecessary in most cases. Instead of using identity $=$ we can employ membership $\in$ in most cases of everyday mathematics. But the fun part of this paper is $Z\Phi$, as we see below.

Observe sentence \ref{1} is not a formula of $Z\Phi$, since $Z\Phi$ is a first-order theory without identity. But sentence \ref{2} is indeed a formula of $Z\Phi$. If we assume differential and integral calculus of real functions is grounded on any transitive model of $Z\Phi$, then, according to Propositions \ref{meta} and \ref{meta2}, sentences \ref{1} and \ref{2} are semantically equivalent (from the perspective of any transitive model of $Z\Phi$, as long such a model is endowed with identity).

It is quite clear that sentence \ref{2} is more difficult to write than sentence \ref{1}. But that's not a problem. If we assume that differential and integral calculus is grounded on any transitive model of $Z\Phi$, we can rewrite sentence \ref{2} as sentence \ref{1}. But in that case the symbol $=$ is just a metalinguistic symbol providing an easy way to write formulas down.

So, the only important difference between a foundational framework like ZF and our proposal is that the latter demands a metatheoretical constraint, namely, a specific class of models of $Z\Phi$ is required to make standard mathematics as it is. Those models are supposed to be transitive.

Obviously, many mathematicians do not want to think deeply about models of any set theory. All they want is a safe place where they can work. ZF provides that, despite the fact that no one knows for sure if ZF is consistent. But $Z\Phi$ provides that as well, as long we do not use $Z\Phi$ alone, but any transitive model of $Z\Phi$. We can still use the symbol $=$ in $Z\Phi$, as long it is just a metalinguistic symbol which is inevitable in any transitive model of $Z\Phi$ endowed with identity.

\section{Some remarks}

The reader could be suspicious about the claims in this paper for some reasons. For starters, Proposition \ref{meta2} demands a model of $Z\Phi$ endowed with identity. So, anyone could ask: is this just a trick where identity is somehow hidden but always present?

Some people \cite{Bueno-14} say it is impossible to get rid of a very general notion of identity which is far more fundamental than a binary predicate letter $=$. Since an identity relationship in a given model does not need to be formal (like it has to be in any first-order theory with identity), we could assume that every model of any theory is endowed with some sort of identity. Nevertheless, the notion of identity in a paper like Bueno's \cite{Bueno-14} is not exactly translatable by an interpretation function in the usual sense of model theory. That happens because the notion of identity in that paper applies even to formulas, something that does not occur in first-order theories. So, let's stick to the usual sense of identity where reflexivity, symmetry, transitivity, and substitutivity between {\em terms\/}, and only terms, are supposed to hold. In that case, yes, identity is still there, since Proposition \ref{meta2} demands the use of transitive models of $Z\Phi$ which are endowed with identity. But, first of all, if we write something like $x=y$, that is not a formula of $Z\Phi$, since identity symbol $=$ belongs to a transitive model of $Z\Phi$. Such a symbol does not belong to the vocabulary of $Z\Phi$. Second, a formula like $\forall t(t\in x \Leftrightarrow t\in y)$ {\em is\/} a formula of $Z\Phi$. And such a formula is equivalent to $x=y$ in any transitive model of $Z\Phi$ endowed with identity. Henceforth, we have at least two possible ways to cope with $Z\Phi$: 

\begin{itemize}

\item[(i)] we write standard mathematics formulas using only formulas of $Z\Phi$ and, besides that, we assume that those formulas refer to transitive models of $Z\Phi$ with identity; or

\item[(ii)] we write standard mathematics formulas using only translations of formulas of $Z\Phi$ into formulas of a transitive model of $Z\Phi$ with identity. 

\end{itemize}

In the first case it is impossible to write equations. In the second case it is possible. But in both cases equations are unnecessary. And in both cases we have a foundation for standard mathematics. Nevertheless, option (i) above is somehow confusing, since it does not let clear how this relation between $Z\Phi$ and a transitive model of $Z\Phi$ works. So, a better alternative is option (ii): standard mathematics is done within a transitive model of $Z\Phi$.

On the other hand, and that's very important, there is a third way to cope with $Z\Phi$, which is a subtle variation of option (ii) above. The very concept of transitivity of any model ${\cal U} = (U, I)$ of $Z\Phi$ {\em does not require\/} identity $=$. After all, ${\cal U} = (U, I)$ is transitive iff for any $x$ and any $y$, `$x$ belonging to $U$ and $y$ belonging to $x$' entail `$y$ belonging to $U$', where `belonging' is the translation of $\in$ by means of function $I$. Thus, transitivity, by itself, is enough for all practical purposes. In this sense, {\em identity does not need to be mentioned, ever\/}.

Another possible criticism is the dependence upon metamathematical considerations in order to use $Z\Phi$ as a foundation for standard mathematics. Well, in a way, that's not so unusual. When we admit ZF or variations as a foundation for mathematics we always assume a metamathematical assumption, namely, the consistency of ZF. That takes place because there is no proof of either consistency or inconsistency of ZF. The metamathematical assumption of consistency of ZF is far more radical than my proposal, since it is based on our ignorance about ZF's consistency.

On the other hand, an easier way to guarantee our metamathematical assumption concerning transitive models of $Z\Phi$ is by assuming identity and the Axiom of Extensionality, a maneuver far more simple to adopt. That's the way of ZF. But what is the real advantage of using identity $=$ and Extensionality? The answer is this: to grant an equivalence relation where substitutivity holds! Nevertheless, in a pragmatical sense, that is accomplished with biconditional operator as well: if $F$, $G$, and $H$ are formulas, then $F\Leftrightarrow F$ is a theorem (reflexivity) ; $F\Leftrightarrow G$ entails $G\Leftrightarrow F$ (symmetry); and the conjunction of both $F\Leftrightarrow G$ and $G\Leftrightarrow H$ entails $F\Leftrightarrow H$ (transitivity). Besides, if $F\Leftrightarrow G$ is a theorem and $G$ occurs in any formula, we can replace it by $F$ in the same formula (substitutivity). The obvious problem with that analogy is that it assumes quantification over formulas, something that does not happen in first-order theories. Besides, formulas are supposed to say something about terms. And terms are the objects of mathematical study and investigation. That is why I propose the use of biconditionals to say something about terms $x$ and $y$, like formula $\forall t(t\in x\Leftrightarrow t\in y)$ instead of $x = y$. Naturally I am not suggesting we should use membership formulas instead of equations. I am simply showing that we can do standard mathematics without any equation at all.

\section{Quantum non-individuality}

Some physical phenomena and physical theories have inspired the production of a vast philosophical literature about metaphysical implications of quantum mechanics. One of the issues raised by philosophers and even physicists is the so-called {\em problem of identity and individuality\/}. For details, see, e.g., \cite{French-06} and the bibliography cited in that book. That problem found its place in philosophical discussions due to two factors: (i) if we accept the existence of quantum particles (in some sense of the word `particle'), there is experimental evidence on the existence of multiple indiscernible (in a strong sense of the word) particles; (ii) within the scope of usual formalisms of quantum mechanics, it is necessary to assume some sort of strong notion of indiscernibility among quantum particles; otherwise we cannot handle concepts like quantum distributions and quantum entanglement, among others. Since those two factors above seem to violate Leibniz’s famous Principle of the Identity of Indiscernibles, then any metaphysical inquiry about individuality in the quantum realm is simply inevitable. Everyday objects like chairs and apples do not exhibit such a radical notion of indiscernibility. If apple $a$ is indiscernible from apple $b$, with respect to all features there is to know about apples, then we are talking about one single apple who admits two names: $a$ and $b$. On the other hand, there may exist two or more truly indiscernible electrons. In a nutshell, we can label apples, but not electrons.

Some philosophers, for example, assume that quantum particles cannot be regarded as individuals. Within this context some non-standard formal theories have been proposed in order to provide a detailed and precise account of concepts like identity, individuality, and indiscernibility. The partial aftermath of that is that now we have philosophical discussions about such non-standard formal systems as well. So, whoever enters this world of philosophy of physics and philosophy of mathematics, beware! It's not easy.

In this Section I focus on one of the non-standard formal systems proposed, namely, quasi-set theory \cite{Krause-92}. Although there are different versions of quasi-set theory \cite{Krause-92}\cite{Krause-05}\cite{French-06}, that's not an issue for my purposes. That's because all those versions share some basic notions and axioms in common. One of the most important concepts in quasi-set theories is that one of indiscernibility (indistinguishability), usually denoted by $\equiv$. If $a$ and $b$ are terms, then we read $a\equiv b$ as `$a$ and $b$ are indiscernible (indistinguishable)'. The idea is that $\equiv$ is an equivalence relation where substitutivity does not necessarily hold. Thus, $\equiv$ is supposed to be weaker than identity $=$.

Well, in the formal framework proposed here, $Z\Phi$, it is possible to introduce an equivalence relation where substitutivity does not necessarily hold as well. We can call it {\em similarity\/} and denote it by $\cong$, as it follows:

\begin{definition}\label{defineconga}

$a\cong b$ iff $\forall t(t\in a \Leftrightarrow t\in b)$.

\end{definition}

In other words, $a$ and $b$ are {\em similar\/} iff they share the same elements. Observe that this informal notion of `same elements' has nothing to do with identity in the usual sense of the predicate letter $=$ in first-order theories with identity. The notion of `same elements' is informally captured by {\em two\/} occurrences of term $t$ in the formula $t\in a \Leftrightarrow t\in b$. So, the assessment of {\em two\/} occurrences of the {\em same\/} term in the {\em same\/} formula is based on a metalinguistic consideration. So, my informal reading (in English) of Definition \ref{defineconga} is not the issue here. The point here is the {\em definiens\/} (a formula!) in Definition \ref{defineconga}. That's very important to be crystal clear. After all, there are different concepts of identity in the literature of philosophy. Bueno's work \cite{Bueno-14} is a great example to illustrate this diversity of ideas about identity. In that paper there is a very radical notion of identity applicable even to formulas and auxiliary symbols like comma and parenthesis. Here, nevertheless, identity refers solely to a predicate letter in first-order theories with identity. English words like `same', `identical', `equal', `self', `alike' and other synonyms are no translation of a formal concept like the binary predicate $=$. Sound translations of primitive concepts of a formal theory are provided by model theory. Quite frequently people get confused when they try to translate a formal language into something like English.

A similar situation occurs in quasi-set theory. Since indistinguishability $\equiv$ is an equivalence relation, the first postulate of that theory says this: $\forall x(x\equiv x)$. Informally speaking, every term is indistinguishable from {\em itself\/}. 

It's easy to show that our binary relation $\cong$ from last definition is an equivalence relation. If ${\cal U} = (U, I)$ is a non-transitive model of $Z\Phi$, then the translation $I(\cong)$ is an equivalence relation in $U$ where substitutivity does not necessarily hold. If ${\cal U} = (U, I)$ is a transitive model of $Z\Phi$, then the translation $I(\cong)$ is an equivalence relation in $U$ where substitutivity holds. In that case, $I(\cong )$ is simply the usual notion of identity in $U$. One obvious advantage of our semantical approach is that it demands one single primitive concept, namely, membership $\in$. In quasi-set theory indistinguishability is not defined. It is assumed as a primitive concept, together with membership $\in$. Besides, in quasi-set theory there are other primitive concepts like quasi-cardinality, micro-atoms, macro-atoms, a unary predicate letter $Z$ (to discriminate the so-called sets), and a huge list of 28 axioms in one of its versions \cite{French-06}. Besides, quasi-set theory usually employs the notion of atoms. That entails the existence of standard permutation models with identity \cite{Sant'Anna-23}.

So, how to use $Z\Phi$ to deal with ensembles of multiple indistinguishable (in a strong sense of violation of LPII) quantum particles and still keep the standard mathematics of Hilbert spaces and probability theory so frequently used to describe quantum mechanics? One obvious possibility is the use of some specific cases of non-transitive models of $Z\Phi$. How to either select or to build those models? Here follows a recipe.

\begin{enumerate}

\item Let ZFA be Zermelo-Fraenkel set theory with a finite set of atoms \cite{Jech-03}. Remember ZFA has a monadic predicate `to be a set', and an atom is an empty term who is not a set. Thus, Extensionality does not apply to atoms in ZFA. 
    
\item Let ${\cal U} = (s\cup a, I)$ be a model of ZFA, where $s$ is a proper class of sets and $a$ is the corresponding finite set of atoms in ${\cal U}$, in the sense that no object from $s\cup a$ belongs to any member of $a$ (within the model). Besides, the transitive closure of $\in$ on any term from $s$ has no element from $a$. For all that matters ${\cal U}$ could be a very standard and boring permutation model of ZFA \cite{Jech-03}. Boring in the sense of a permutation group whose $\in$-automorphisms are all trivial. 

\item If ${\cal U}$ is, for example, a permutation model of ZFA like the one sketched in last item above, then there is one single term $\emptyset$ belonging to $s$ such that $\emptyset$ is empty. After all, ${\cal U}$ is a transitive model of ZFA. The difference between $\emptyset$ and any term from $a$ is that $\emptyset$ belongs to $s$ but not to $a$. Besides, no term from $a$ belongs to $s$.

\item If $\in_{\cal U}$ is the translation $I(\in)$ of membership in ZFA, and $=$ is the translation $I(=)$ of identity in ZFA, then we can replace all formulas $x = y$ by $\forall t(t\in_{\cal U} x \Leftrightarrow t\in_{\cal U} y)$. That does not change the fact that ${\cal U}$ is transitive. 
    
\item Thanks to last item we can drop identity $=$ from ${\cal U}$, in the sense that translation function $I$ is restricted to one single relation, namely, membership. Such a restriction of $I$ is denoted by $I'$.

\item Now, let ${\cal U}' = (s\cup \wp(a), I')$ be an interpretation of $Z\Phi$, where $\wp(a)$ is the power of $a$. Clearly ${\cal U}'$ is a non-transitive model of $Z\Phi$, since any subset of $a$ that does not belong to $s$ is empty. For example, $a$ itself is empty since no term in $a$ belongs to the universe $s\cup \wp(a)$.
    
\item On the other hand ${\cal U'}$ admits a transitive sub-model $(s,I')$. 

\item So, when dealing with standard mathematics we use sub-model $(s,I')$. When dealing with ensembles of objects that demand violation of LPII, like quantum particles, we use the non-transitive counterpart of ${\cal U}'$. Objects $x$ and $y$ are indiscernible iff they share the same elements from the universe $s\cup \wp(a)$. That's equivalent to say they do not share the same elements from universe $s\cup \wp(a)$. If $a_1$ is an atom from $a$ and $a_2$ is an atom from $a$, then both $\{a_1\}$ and $\{a_2\}$ belong to $s\cup \wp(a)$, and $\{a_1\}$ and $\{a_2\}$ are indiscernible, since neither $a_1$ nor $a_2$ belongs to $s\cup \wp(a)$. Nevertheless, $\{a_1\}$ and $\{a_2\}$ are not necessarily identical since ${\cal U}'$ is non-transitive. The model itself `knows' that $\{a_1\}$ and $\{a_2\}$ are different if $a_1$ is different from $a_2$. But such a difference is not translatable by function $I'$. 

\end{enumerate}

An analogous recipe could be adapted for dealing with permutation models of quasi-set theory. Among terms of that theory there are two kinds of atoms: micro-atoms and macro-atoms. Macro-atoms are supposed to behave in a similar way to atoms in ZFA, while micro-atoms behave somehow differently. That difference is expressed by means of an indistinguishability relationship $\equiv$ which is an equivalence relation. Indistinguishability collapses to identity when dealing with macro-atoms. But that does not happen when dealing with micro-atoms. Thus, if $a_M\cup a_m$ is the class of all atoms in a permutation model of ZFA, such that $a_M\neq\emptyset$, $a_m\neq\emptyset$, and $a_M\cap a_m = \emptyset$, then we can induce a new model of quasi-set theory where we drop $a_m$ (in a similar way we dropped $a$ in model ${\cal U}'$ above). Nevertheless, that's a more complicated maneuver, since quasi-set theory has more primitive concepts and more axioms. Since it was already proved that quasi-set theory admits the existence of permutation models with identity \cite{Sant'Anna-23}, one natural way to mathematically prove that LPII is actually violated is by means of the employment of non-transitive models of quasi-set theory.

\section{Acknowledgments}

This paper was inspired on a conversation I had with Ot\'avio Bueno regarding his paper \cite{Bueno-14}.





\end{document}